\documentclass[12pt,reqno]{article}
\usepackage{amssymb}
\usepackage{amsmath}

\newcommand{\be}{\begin{equation}}
\newcommand{\ee}{\end{equation}}
\newcommand{\bea}{\begin{eqnarray*}}
\newcommand{\eea}{\end{eqnarray*}}
\newcommand{\ba}{\begin{array}}
\newcommand{\ea}{\end{array}}
\newcommand{\bi}{\begin{itemize}}
\newcommand{\ei}{\end{itemize}}
\newcommand{\bc}{\begin{center}}
\newcommand{\ec}{\end{center}}
\newcommand{\bfr}{\begin{flushright}}
\newcommand{\efr}{\end{flushright}}

\newtheorem{theorem}{Theorem}[section]

\newtheorem{corollary}[theorem]{Corollary}

\newtheorem{definition}[theorem]{Definition}

\numberwithin{equation}{section}

\begin{document}

\title{\bf On the stability of the linear functional equation
           in a single variable on complete metric groups}
\author{Soon-Mo Jung, Dorian Popa, Michael Th. Rassias}
\date{\small Mathematics Section,
             College of Science and Technology, \\
             Hongik University,
             339--701 Sejong, Republic of Korea \\
             E-mail: smjung@hongik.ac.kr \\
             \medskip
             Department of Mathematics,
             Technical University of Cluj-Napoca, \\
             28 Memorandumului, 400114 Cluj-Napoca,
             Romania \\
             E-mail: Popa.Dorian@math.utcluj.ro \\
             \medskip
             Department of Mathematics,
             ETH--Z\"{u}rich, Ramistrasse 101, \\
             8092 Z\"{u}rich, Switzerland \\
             E-mail: michail.rassias@math.ethz.ch}
\maketitle

\begin{quote}
{\bf Abstract.}
In this paper we obtain a result on Hyers-Ulam stability of
the linear functional equation in a single variable
$f(\varphi(x)) = g(x) \cdot f(x)$ on a complete metric group.
\end{quote}

\section{Introduction}

\footnotetext{\textit{Mathematics Subject Classification (2010)}:
              65J15.}
\footnotetext{\textit{Key words and phrases}:
              Nonlinear operator, stability,
              functional equation, complete metric group,
              inequalities, Banach spaces, operator mapping,
              Euler-Mascheroni constant.}

Hyers-Ulam stability is one of the main topics in the theory
of functional equations.
Generally a functional equation is said to be stable provided,
for any function $f$ satisfying the perturbed functional
equation, there exists an exact solution $f_0$ of that equation
which is not far from the given $f$.
Based on this concept, the study of the stability of functional
equations can be regarded as a branch of optimization theory.
(We can find some applications of the Hyers-Ulam stability to
optimization theory and economics in \cite{castillo}.)

It seems that the first result on the stability of functional
equations appeared in the famous book by Gy. P\'{o}lya and G.
Szeg\"{o} \cite{13} and concerns the Cauchy functional equation
on the set of positive integers.
But the starting point of the stability theory of functional
equations is due to S.M. Ulam who formulated a question
concerning the perturbation of homomorphisms on metric groups.
The first result for Ulam's problem was obtained by D.H. Hyers
for the Cauchy functional equation on Banach spaces.
Due to the question of Ulam and the answer of Hyers the stability
of functional equations is called after their names.

For more details on Hyers-Ulam stability of functional equations
and optimization theory we refer the reader to
\cite{1,3,cesaro,czerw,11,pardalos,rocka}.

The functional equation
\begin{equation}
\label{1.1}
f(\varphi(x)) = g(x)f(x) + h(x),
\end{equation}
where $f$ is the unknown function and $g, h, \varphi$ are given
functions, is called the linear functional equation in a single
variable.
For particular cases of $g$ and $h$ in (\ref{1.1}) we obtain
some classical functional equations.
We mention here some of them as
\begin{itemize}
\item Abel's equation
\begin{equation}
\label{1.2}
f(\varphi(x)) = f(x) + c
\end{equation}

\item Schr\"{o}der's equation
\begin{equation}
\label{1.3}
f(\varphi(x)) = cf(x)
\end{equation}

\item Gamma functional equation
\begin{equation}
\label{1.4}
f(x+1) = xf(x)
\end{equation}

\item Digamma functional equation
\begin{equation}
\label{1.5}
f(x+1) = f(x) + \frac{1}{x}.
\end{equation}
\end{itemize}

Recall that Digamma function
$\psi _0 : \mathbb{R}_+^* \to \mathbb{R}$ is defined by
\begin{equation}
\label{1.6}
\psi_0(x) = \frac{d}{dx} \ln \Gamma(x)
          = \frac{\Gamma'(x)}{\Gamma(x)},
\hspace*{5mm} \forall\ x \in \mathbb{R}_+^*,
\end{equation}
where
\begin{equation}
\label{1.7}
\Gamma(x) = \int_0^\infty t^{x-1} e^{-t} dt,
\hspace*{5mm} \forall\ x \in \mathbb{R}_+^*,
\end{equation}
$\mathbb{R}_+$ stands for the set of all nonnegative numbers,
i.e., $\mathbb{R}_+ = [0, \infty)$ and
$\mathbb{R}_+^* = (0, \infty)$.
For more details on the functional equation (\ref{1.1}) and
its particular cases we refer to \cite{12} and the references
therein.
It seems that the first result on stability for the equation
(\ref{1.1}) was obtained in 1970 by J. Brydak \cite{2}.
A generalized Hyers-Ulam stability of the gamma functional
equation was obtained by S.-M. Jung in \cite{10}.
A nice result on generalized Hyers-Ulam stability of the
equation (\ref{1.1}) was obtained by T. Trif \cite{14} for
functions $f$ acting from an arbitrary nonempty set $S$ into a
Banach space $X$.

Some recent results on the stability and nonstability of the
equation (\ref{1.1}) and the linear functional equation of
higher order in a single variable were obtained by J. Brzdek,
D. Popa, B. Xu (see \cite{5,6,7,8,9}).

The goal of this paper is to study the Hyers-Ulam stability of the
homogeneous linear functional equation (\ref{1.1}) for functions
defined from an arbitrary nonempty set $S$ into a complete metric
group $(G, \cdot, d)$, i.e., $(G, \cdot)$ is a group, $(G, d)$ is
a complete metric space, the group's binary  operation and the
inverse operation are continuous with respect to the product
topology on $G \times G$ and the topology generated by the metric
$d$ on $G$, respectively.

\section{Stability of linear functional equation}

Let $S$ be a nonempty set, $(G, \cdot, d)$ a complete metric group
with the metric $d$ invariant to left translations, i.e.,
\begin{equation}
\label{2.1} d(x \cdot y, x \cdot z)= d(y, z), \hspace*{5mm}
\forall\ x, y, z \in G,
\end{equation}
and let $\varphi : S \to S$, $g : S \to G$ be given functions.
An example of metric invariant to left translations is the
metric induced by a norm.

We deal with the Hyers-Ulam stability of the linear functional
equation
\begin{equation}
\label{2.2}
f(\varphi(x)) = g(x) \cdot f(x),
\end{equation}
where $f : S \to G$ is the unknown function.

Let $\mathbb{R}_+^S$ be the class of all functions $\varepsilon :
S \to \mathbb{R}_+$. We study the generalized Hyers-Ulam stability
of the equation (\ref{2.2}) in the sense defined in \cite{3}.

\begin{definition}
{\rm Let ${\mathcal C} \subseteq \mathbb{R}_+^S$ be nonempty
and ${\mathcal T}$ be an operator mapping ${\mathcal C}$ into
$\mathbb{R}_+^S$.
We say that the equation (\ref{2.2}) is ${\mathcal T}$-stable
(with uniqueness, respectively) provided for every
$\varepsilon \in {\mathcal C}$ and $f : S \to G$ with
$$
d(f(\varphi(x)), g(x) \cdot f(x)) \leq \varepsilon(x),
\hspace*{5mm} \forall \ x \in S
$$
there exists a (unique, respectively) solution $f_0 : S \to G$
of the equation (\ref{2.2}) such that
$$
d(f(x), f_0(x)) \leq {\mathcal T} \varepsilon(x),
\hspace*{5mm} \forall \ x \in S.
$$
}
\end{definition}
If $\varepsilon$ is a constant function in the previous definition
then the equation (2.2) is said to be stable in Hyers-Ulam sense.

By $\varphi^k$, $k \in \mathbb{N}_0 = \mathbb{N} \cup \{ 0 \}$
we denote the $k$-th iterate of the function $\varphi$,
$\varphi^0 = 1_S$, $\varphi^k = \varphi \circ \varphi^{k-1}$,
$k \in \mathbb{N}$.

The main result is contained in the next theorem.

\begin{theorem}
\label{t2.2}
Let $\varepsilon : S \to \mathbb{R}_+$ be a function with the
property
\begin{equation}
\label{2.3}
\sum_{n=0}^\infty \varepsilon(\varphi^n(x)) = \Phi(x),
\hspace*{5mm} \forall \ x \in S,
\end{equation}
where $\Phi : S \to \mathbb{R}_+$.
Then for every function $f : S \to G$ satisfying the inequality
\begin{equation}
\label{2.4}
d(f(\varphi(x)), g(x) \cdot f(x)) \leq \varepsilon(x),
\hspace*{5mm} \forall \ x \in S,
\end{equation}
there exists a unique solution $f_0 : S \to G$ of the functional
equation (\ref{2.2}) such that
\begin{equation}
\label{2.5}
d(f(x), f_0(x)) \leq \Phi(x), \hspace*{5mm} \forall \ x \in S.
\end{equation}
\end{theorem}

\noindent
{\it Proof.}
{\it Existence.}
Let $f : S \to G$ be a function satisfying (\ref{2.4}).
Then the following relation holds:
\begin{equation}
\label{2.6}
d\!\left( f\big( \varphi^n(x) \big),
          \prod_{k=1}^n g\big( \varphi^{k-1}(x) \big) \cdot f(x)
   \right)
\leq \sum_{k=1}^n \varepsilon\big( \varphi^{k-1}(x) \big)
\end{equation}
for all $x \in S$ and $n \in \mathbb{N}$.
We prove (\ref{2.6}) by induction on $n$.
Since the group $(G, \cdot)$ is not generally commutative, we
let
$$
\prod_{k=p}^n a_k := a_n \cdot a_{n-1} \cdot \ldots \cdot a_p,
$$
where $a_k \in G$ for $p \leq k \leq n$.

For $n=1$ the relation (\ref{2.6}) holds in view of (\ref{2.4}).
We suppose that (\ref{2.6}) holds for some $n \in \mathbb{N}$
and for all $x \in S$, and we prove that
$$
d\!\left( f\big( \varphi^{n+1}(x) \big),
          \prod_{k=1}^{n+1} g\big( \varphi^{k-1}(x) \big) \cdot
          f(x)
   \right)
\leq \sum_{k=1}^{n+1} \varepsilon\big( \varphi^{k-1}(x) \big),
\hspace*{5mm} x \in S.
$$
Indeed, it follows from (\ref{2.1}), (\ref{2.4}) and (\ref{2.6})
that
\begin{align*}
d&\!\left( f\big( \varphi^{n+1}(x) \big),
           \prod_{k=1}^{n+1} g\big( \varphi^{k-1}(x) \big)
           \cdot f(x)
    \right) \\
&\leq d\big( f\big( \varphi^{n+1}(x) \big),
             g(\varphi^n(x)) \cdot f(\varphi^n(x)) \big) \\
&\hspace*{5mm}
      + d\!\left( g(\varphi^n(x)) \cdot f(\varphi^n(x)),
                  \prod_{k=1}^{n+1}
                  g\big( \varphi^{k-1}(x) \big) \cdot f(x)
           \right) \\
&\leq \varepsilon(\varphi^n(x)) +
      d\!\left( f(\varphi^n(x)),
                \prod_{k=1}^n g\big( \varphi^{k-1}(x) \big)
                \cdot f(x)
         \right) \\
&\leq \sum_{k=1}^{n+1} \varepsilon\big( \varphi^{k-1}(x) \big),
      \hspace*{5mm} x \in S.
\end{align*}
Hence (\ref{2.6}) holds for all $x \in S$ and $n \in \mathbb{N}$.

Now let $( \varepsilon_n )_{n \geq 1}$ be the sequence of
functions defined by
\begin{equation}
\label{2.7}
\varepsilon_n(x) :=
\left( \prod_{k=1}^n g\big( \varphi^{k-1}(x) \big) \right)^{-1}
\cdot f(\varphi^n(x)),
\hspace*{5mm} n \in \mathbb{N},\ x \in S.
\end{equation}
We prove that $( \varepsilon_n(x) )_{n \geq 1}$ is a Cauchy
sequence in $(G, \cdot, d)$ for all $x \in S$, where $a^{-1}$
means the inverse of the element $a$ in the group $G$.
Using (\ref{2.1}) and (\ref{2.6}), we have
\begin{align}
d&\big( \varepsilon_{n+p}(x), \varepsilon_n(x) \big) \nonumber\\
&= d\!\left(\! \left( \prod_{k=1}^{n+p}
                      g\big( \varphi^{k-1}(x) \big)\!
               \right)^{-1}\! \cdot f\big( \varphi^{n+p}(x) \big),
               \left( \prod_{k=1}^n g\big( \varphi^{k-1}(x) \big)\!
               \right)^{-1}\! \cdot f(\varphi^n(x))
    \!\right) \nonumber\\
&= d\!\left(\! \left( \prod_{k=n+1}^{n+p}
                         g\big( \varphi^{k-1}(x) \big)\!
                  \right)^{-1}\! \cdot
                  f\big( \varphi^{n+p}(x) \big),
                  f(\varphi^n(x))
       \!\right) \nonumber\\
&\leq \sum_{k=1}^{p}
      \varepsilon\big( \varphi^{k-1}(\varphi^n(x)) \big)
 \leq \sum_{k=0}^\infty \varepsilon\big( \varphi^{n+k}(x) \big)
      \label{2.8}
\end{align}
for $x \in S$ and $n, p \in \mathbb{N}$.

Now
$r_n(x) := \sum_{k=0}^\infty \varepsilon(\varphi^{n+k}(x))$,
$n \in \mathbb{N}$, is the remainder of order $n$ of the
convergent series (\ref{2.3}), so $\lim_{n \to \infty} r_n(x) = 0$
for all $x \in S$.
We conclude that $( \varepsilon_n(x) )_{n \geq 1}$ is a Cauchy
sequence, therefore it is convergent since $G$ is a complete
metric group.
Define the function $f_0$ by
$$
f_0(x) = \lim_{n \to \infty} \varepsilon_n(x),
\hspace*{5mm} x \in S.
$$
The relation (\ref{2.8}), for $p=1$, leads to
\begin{equation}
\label{2.9} d( \varepsilon_{n+1}(x), \varepsilon_n(x) ) \leq
\sum_{k=0}^\infty \varepsilon\big( \varphi^{n+k}(x) \big),
\hspace*{5mm} n \in \mathbb{N},\ x \in S.
\end{equation}

Taking account of $\varepsilon_{n+1}(x) = g(x)^{-1} \cdot
\varepsilon_n(\varphi(x))$ and letting $n \to \infty$ in
(\ref{2.9}) it follows that
$$
d\big( g(x)^{-1} \cdot f_0(\varphi(x)), f_0(x) \big) = 0
$$
which is equivalent to
$f_0(\varphi(x)) = g(x) \cdot f_0(x)$, $x \in S$, i.e., $f_0$
is a solution of the equation (\ref{2.2}).

On the other hand, the relations (\ref{2.1}) and (\ref{2.6})
lead to
\begin{equation}
\label{2.10}
d(\varepsilon_n(x), f(x))
\leq \sum_{k=1}^n \varepsilon\big( \varphi^{k-1}(x) \big)
\end{equation}
for all $x \in S$ and $n \in \mathbb{N}$, therefore letting
$n \to \infty$ in (\ref{2.10}) we get
$$
d(f_0(x), f(x)) \leq \Phi(x),
$$
which completes the proof of the existence.

{\it Uniqueness.} Assume that for a function $f$ satisfying
(\ref{2.4}) there exist two solutions $f_1, f_2$ of the equation
(\ref{2.2}) satisfying
$$
d(f(x), f_i(x)) \leq \Phi(x),
\hspace*{5mm} \forall\ x \in S,\ i \in \{ 1, 2 \}
$$
and $f_1 \neq f_2$.

Taking into account that $f_1, f_2$ satisfy (\ref{2.2}), it
follows easily that
$$
f_i(\varphi^n(x))
= \prod_{k=1}^n g\big( \varphi^{k-1}(x) \big) \cdot f_i(x),
\hspace*{5mm} n \in \mathbb{N},\ x \in S,\ i \in \{ 1, 2 \},
$$
and hence
\begin{align*}
d&(f_1(x), f_2(x)) \\
&= d\!\left(\! \left( \prod_{k=1}^n
                      g\big( \varphi^{k-1}(x) \big)\!
               \right)^{-1}\! \cdot f_1(\varphi^n(x)),
               \left( \prod_{k=1}^n
                      g\big( \varphi^{k-1}(x) \big)\!
               \right)^{-1}\! \cdot f_2(\varphi^n(x))
    \!\right) \\
& = d\big( f_1(\varphi^n(x)), f_2(\varphi^n(x)) \big) \\
&\leq d\big( f_1(\varphi^n(x)), f(\varphi^n(x)) \big) +
      d\big( f(\varphi^n(x)), f_2(\varphi^n(x)) \big) \\
&\leq 2\Phi(\varphi^n(x)),
      \hspace*{5mm} x \in S,\ n \in \mathbb{N}.
\end{align*}
Since
$\lim_{n \to \infty} \Phi(\varphi^n(x)) = \lim_{n \to \infty}
 r_n(x) = 0$, $x \in S$, it follows that $f_1(x) = f_2(x)$,
which completes the proof.
\hfill$\Box$
\vspace{5mm}

The Digamma function $\psi_0 : \mathbb{R}_+^* \to \mathbb{R}$
is defined by (\ref{1.6}).
The Digamma function is frequently called the psi function and
it satisfies the Digamma functional equation (\ref{1.5}) for
all $x \in \mathbb{R}_+^*$.
Indeed, we know that $\psi_0$ is the unique solution of the
functional equation (\ref{1.5}) which is monotone on
$\mathbb{R}_+^*$ and satisfies $\psi_0(1) = -\gamma$, where
$\gamma = 0.577215\ldots$ is the Euler-Mascheroni constant
(see \cite[\S 6.3]{abram} and \cite[\S 6.11.5]{zwillinger}).

The gamma function defined by (\ref{1.7}) satisfies the
functional equation
$$
\Gamma(x+1) = x \Gamma(x),
\hspace*{5mm} \forall\ x \in \mathbb{R}_+^*.
$$
If we take the logarithmic values from both sides of the last
equation, then we have
$$
\ln \Gamma(x+1) = \ln \Gamma(x) + \ln x,
\hspace*{5mm} \forall\ x \in \mathbb{R}_+^*.
$$
We differentiate each side of the above equality with respect
to $x$ to get
$$
\frac{d}{dx} \ln \Gamma(x+1)
= \frac{d}{dx} \ln \Gamma(x) + \frac{1}{x},
\hspace*{5mm} \forall\ x \in \mathbb{R}_+^*.
$$
In view of (\ref{1.6}), we know that the Digamma function
$\psi_0$ is a solution of the Digamma functional equation
(\ref{1.5}).

The generalized Hyers-Ulam stability of the Digamma functional
equation (\ref{1.5}) follows from Theorem \ref{t2.2}.

\begin{corollary}
\label{c2.3}
Let $\varepsilon : \mathbb{R}_+^* \to \mathbb{R}_+$ be a
function with the property
$$
\sum_{n=0}^\infty \varepsilon(x+n) = \Phi(x),
\hspace*{5mm} \forall\ x \in \mathbb{R}_+^*.
$$
Then for every function $f : \mathbb{R}_+^* \to \mathbb{R}$
satisfying
$$
\left| f(x+1) - f(x) - \frac{1}{x} \right| \leq \varepsilon(x),
\hspace*{5mm} x \in \mathbb{R}_+^*
$$
there exists a unique solution
$f_0 : \mathbb{R}_+^* \to \mathbb{R}$ of the equation
(\ref{1.5}) such that
$$
| f(x) - f_0(x) | \leq \Phi(x),
\hspace*{5mm} \forall\ x \in \mathbb{R}_+^*.
$$
\end{corollary}

\noindent
{\it Proof.}
Take $S = \mathbb{R}_+^*$, $\varphi(x) = x + 1$, $G = \mathbb{R}$
with the usual addition and $d$ the Euclidian metric on
$\mathbb{R}$ and $g(x) = 1/x$, $x \in \mathbb{R}_+^*$.
Then the result follows in view of Theorem \ref{t2.2}.
\hfill$\Box$
\vspace{5mm}

\end{document}